\documentclass[a4paper,11pt]{amsart}
\usepackage{hyperref,latexsym}
\usepackage{enumerate}
\usepackage{graphicx}
\usepackage{color}

\theoremstyle{plain}

\theoremstyle{definition}

\theoremstyle{remark}

\begin{document}

\title[]
      {Numerical non-integrability of Hexagonal string billiard}

\date{}
\author{Misha Bialy}
\address{School of Mathematical Sciences, Raymond and Beverly Sackler Faculty of Exact Sciences, Tel Aviv University,
Israel} 
\email{bialy@tauex.tau.ac.il}
\thanks{MB and BY were partially supported by ISF grant 580/20  and DFG 2565/7-1 within the Middle East Collaboration Program.}

\author{Baruch Youssin}
\address{School of Mathematical Sciences, Raymond and Beverly Sackler Faculty of Exact Sciences, Tel Aviv University,
	Israel}
\email{baruch@byoussin.com}


\begin{abstract}
We consider a remarkable $C^2$-smooth billiard table introduced by Hans L.~Fetter \cite{fetter}. It is obtained by the string construction from a regular hexagon for a special value of the length of the string. It was suggested in
\cite{fetter} as a possible counter-example to the Birkhoff-Poritsky conjecture. In this paper, we investigate numerically the behavior of this billiard and find chaotic regions near hyperbolic periodic orbits. They are very small since the billiard table is nearly circular. 
\end{abstract}

\maketitle

\section{Introduction}
 Hans L.~Fetter \cite{fetter} proposed a hexagonal string billiard
 as a possible counterexample to the famous Birkhoff-Poritsky conjecture \cite{Poritsky}, which gained much progress recently \cite{K-S1}\cite{BM}. Fetter \cite{fetter} performed numerical simulations of 70 orbits through 480 iterations for the proposed billiard, showing that these orbits appear to lie on invariant curves.

The purpose of this report is twofold. 

First, we calculate 37 orbits spread out over the upper half of the phase cylinder. We calculate these orbits through at least 200,000 iterations each, and most with 500,000 or more; for the way we estimate the precision of our calculations, see \ref{error.estimation} below. We have discovered that these orbits indeed belong to invariant curves, in agreement with the conclusions of \cite{fetter}. 

In our second goal, we consider hyperbolic periodic orbits, one of the period 12 with the rotation number $\frac{5}{12}$ and the other of the period 27 with the rotation number $\frac{11}{27}$. We find tiny chaotic regions close to these orbits; this is a well-known scenario for the creation of chaos; see \cite{Simo} for an example. Thus, our simulation shows that the hexagonal string billiard is not integrable, contrary to the suggestion of \cite{fetter}. 
Note, however, that the chaotic regions we have found, are very small
($10^{-4}\times10^{-5}$ for the first hyperbolic orbit and $3\cdot 10^{-7}\times10^{-8}$ for the second one) comparing to the billiard size, which is of order 1.

The hexagonal string billiard is constructed by the string construction for the regular hexagon in the plane for a special value of the length $l$ of the string. Generally, the curve $\Gamma$, obtained by the string construction from any convex closed curve $\gamma$ for any length of the string $l>length\ \gamma$, is of class $C^1$. Moreover, the curve $\gamma$ is a caustic with respect to the billiard in $\Gamma$. If $\gamma $ is a regular hexagon then there is a special value of $l$ for which $\Gamma$ is of class $C^2$; in this case, $\Gamma$ consists of 6 pieces of congruent ellipses joined together in a $C^2$ way.
{Remarkably, the curve $\Gamma$ looks very close to a circle; see Fig.~\ref{fig.comparison.with.circle} below. This fact explains that the chaotic regions which we found, are so small.}

Another class of $C^2$-curves with non-smooth convex caustics $\gamma$ is suggested in \cite{AB}. 

Note that recently L. Bunimovich introduced an interesting class of billiards, called flower billiards, which contains the hexagonal billiard. 
\section{Construction of the billiard}
Our approach and the algorithm is based on the symplectic coordinates $(\phi, p)$ in the space of oriented lines in the plane, where $\phi$ is the angle formed by the right normal to the line with the x-coordinate axis, and $p$ is the signed distance from the origin to the oriented line. 
With this approach it is convenient to define the billiard table with the help of the support function $h(\theta)$, 
$$
h(\theta):=\max_{(x,y)\in\Gamma}\{x\cos\theta+y\sin\theta\}.
$$

\begin{figure}[h]
	\centering
	\includegraphics[width=0.5\textwidth]{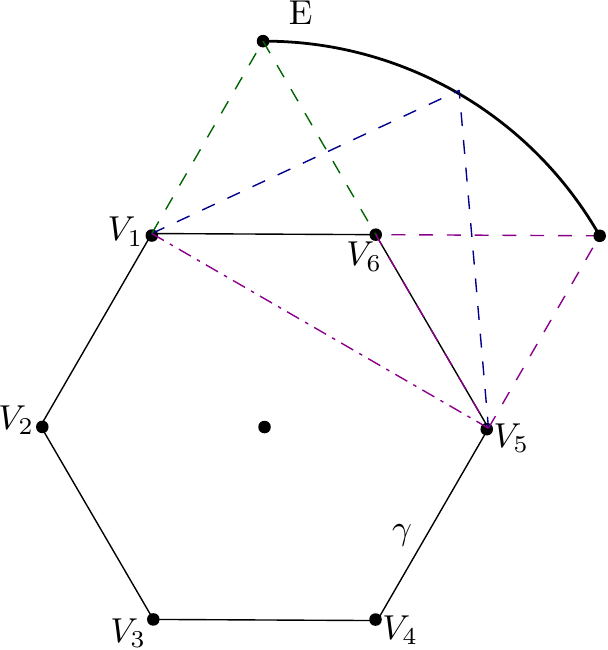}
	\caption{The string construction for $\gamma$ produces 6 congruent elliptic arcs}
	\label{fig.string.construction}
\end{figure}

The hexagonal billiard $\Gamma$ is constructed by the string construction applied to the regular hexagon $\gamma$ with the length of the edge 1, centered at the origin. The length of the string is chosen to be 7; see Fig.~\ref{fig.string.construction} . This choice ensures the $C^2$-smoothness of $\Gamma$, which is the conjunction of six congruent elliptic arcs, whose foci segments are short diagonals of the hexagon $\gamma$.
All these arcs are obtained from one arc $E$ by rotations by $k\pi/3, k=1,..,5.$ Here $E$ is the arc of the ellipse with the foci $V_1,V_5$; this ellipse has half-axes $a=\frac{3}{2}, b=\sqrt{\frac{ 3}{2}}$ and $c=\sqrt{a^2-b^2}=\frac{\sqrt 3}{2}$.
With these parameters, one can easily write the support function $h(\theta)$ of the arc $E$ and then of the whole curve
$\Gamma$ by $\frac{\pi}{3}$-periodicity:
\begin{equation}
\begin {cases}
 h(\theta)=\sqrt{\frac{9}{4}\cos^2(\theta+\frac{\pi}{6})+\frac{3}{2}\sin^2(\theta+\frac{\pi}{6})}+\frac{1}{2}\sin(\theta+\frac{\pi}{6})
 ;\ \theta\in[\frac{\pi}{6},\frac{\pi}{2}]\\

h(\theta+\frac{\pi}{3}):=h(\theta).
\end {cases}
\end{equation}
	\begin{figure}[h]
	\centering
	\includegraphics[width=0.8\textwidth]{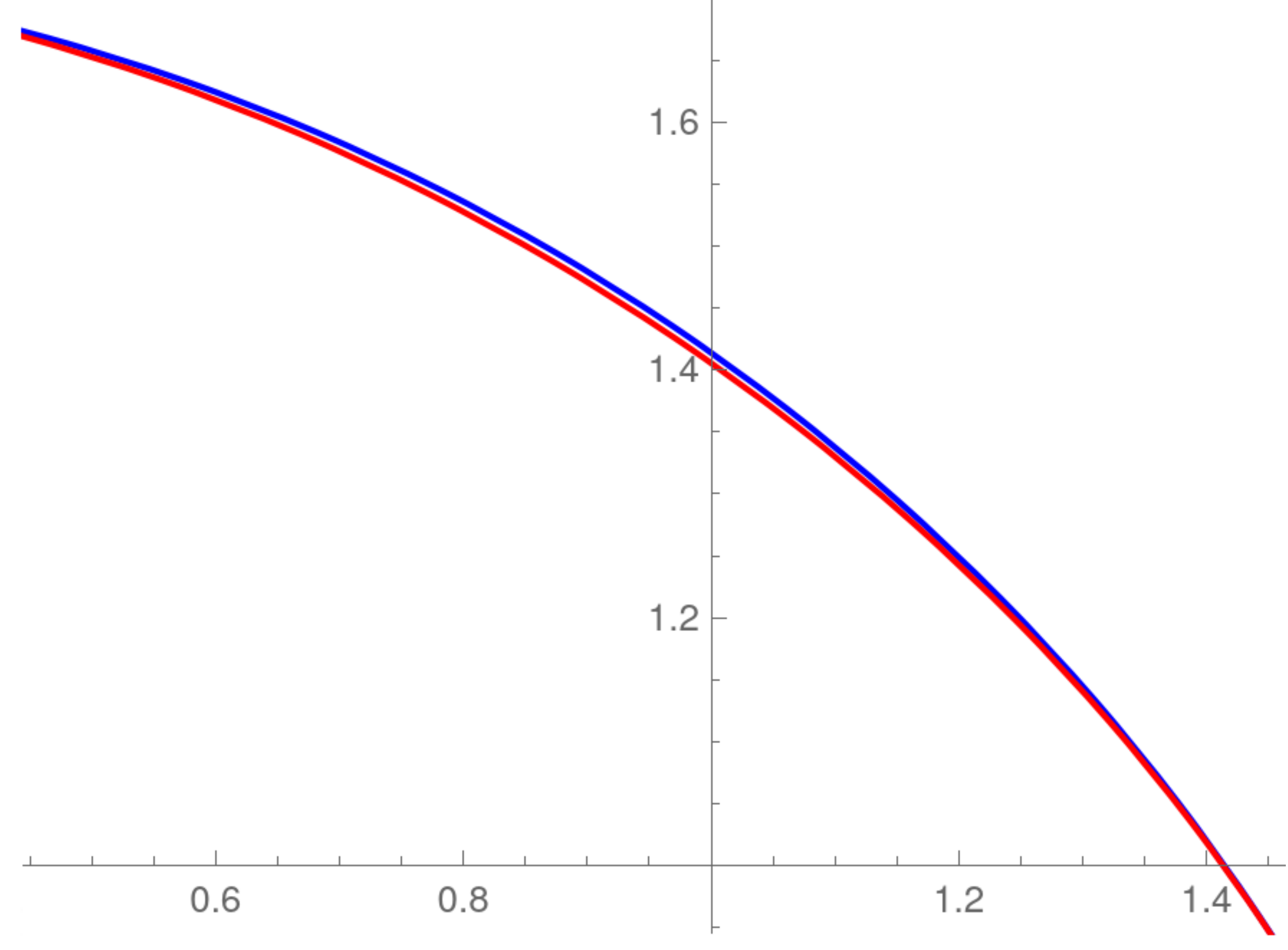}
	\caption{One elliptic arc ({\color{red}\bf red}) compared with the circular arc centered at the origin ({\color{blue}\bf blue})}
	\label{fig.comparison.with.circle}
\end{figure}

\section{The algorithm}

\subsection{Parametric equations of the billiard boundary}
It is well known that the billiard boundary $\Gamma$ can be given by the following parametric equations:
\begin{equation}\label{eq.param.from.support}
\begin{pmatrix}
x\\y
\end{pmatrix}
=h(\theta)
\begin{pmatrix}
\cos\theta\\ \sin\theta
\end{pmatrix}
+h'(\theta)
\begin{pmatrix}
-\sin\theta\\ \cos\theta
\end{pmatrix}
,\ 0\le\theta\le2\pi.
\end{equation}

\subsection{One iteration of the billiard} Consider the action of the billiard map $T$ on the space of oriented lines intersecting $\Gamma$; the coordinates on this space are $(\phi, p)$ as above. Denote $T(\phi, p)$ by $(\phi_1, p_1)$; it can be computed in terms of $(\phi, p)$ in three steps:

	\begin{figure}[h]
	\centering
	\includegraphics[width=0.5\textwidth]{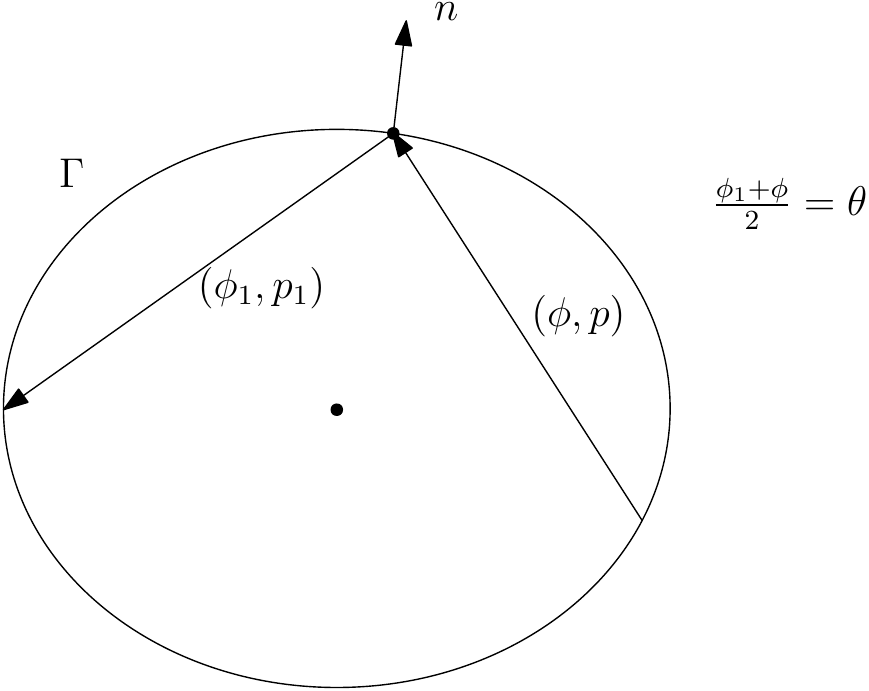}
	\caption{The line $(\phi, p)$ is reflected to $(\phi_1, p_1)$}
	\label{fig.reflection}
\end{figure}

Step 1: Find the parameter value $\theta$ of the \emph{second} of the two intersections of the oriented line $(\phi,p)$ with the billiard boundary.
We find $\theta$ as the solution of the equation
\begin{equation}
x\cos\phi + y \sin\phi = p,
\end{equation}
where $x,y$ are functions of $\theta$ given by the equation (\ref{eq.param.from.support}); to get the \emph{second} intersection (see Fig.~\ref{fig.reflection}), we choose the solution $\theta$ that satisfies $\phi<\theta<\phi+\pi$.

Step 2: Set $\phi_1=2\theta-\phi$.

Step 3: Set $p_1=x\cos\phi_1 + y \sin\phi_1$, where $x,y$ are given by the equation (\ref{eq.param.from.support}).

\subsection{Error estimation} \label{error.estimation}
We apply multiple iterations of the above algorithm and estimate the error as follows.

We perform the calculations with $n$ decimals; to estimate the error, we perform a similar control calculation with $n'$ decimals, $n'\ge n+10$, and take the difference between the results to be the error estimate.

For non-chaotic orbits, we perform the control calculation only through the first 50,000 iterations even if the orbit contains 3,000,000; from our experience, for non-chaotic orbits the error grows with the number of iterations but not significantly.

For non-chaotic orbits, we take $n=30$, $n'=40$ and find that the error on the first 50,000 iterations is no more than $10^{-18}$.

We have performed the control calculations only on a few of the non-chaotic orbits, those close to the billiard boundary and those close to a hyperbolic periodic orbit. As we have considered these orbits most problematic and their error estimates have been consistent, we have not performed control calculations for other orbits.

For the chaotic orbits, we have done the calculations with much higher precision $n$ and have performed the control calculations on entire orbits; see the detailed report below.
\section{Pictures}

	\begin{figure}[h]
	\centering
	\includegraphics[width=0.9\textwidth]{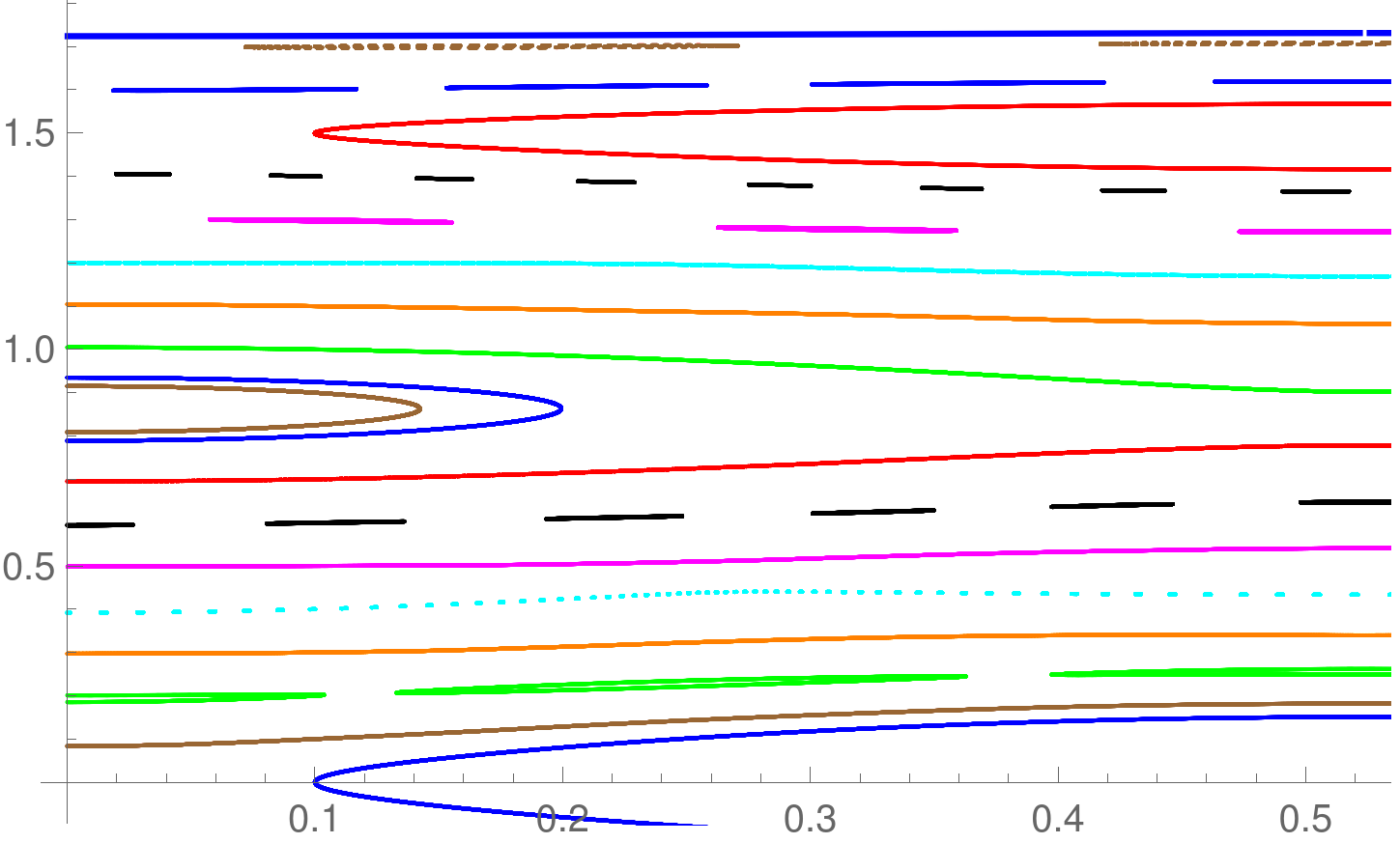}
	\caption{18 orbits spread over the upper half of the phase cylinder. The top blue line is the boundary and the colors are recycled.}
	\label{fig.orbits18spread}
\end{figure}

\subsection{Non-chaotic orbits} \label{subsection.non.chaotic.orbits}
In Fig.~\ref{fig.orbits18spread} we show 18 orbits with the starting points $\phi=0.1$, $p=0,0.1,0.2,...,1.7$; they are spread over the upper half of the phase cylinder.  These orbits were calculated through 500,000 iterations with 30 decimals precision.  We plot only the region $[0, \pi/6]$ on the $\phi$ axis (it is the horizontal axis on all our plots) as the rest is obtained by the symmetry of the billiard.

	\begin{figure}[h]
	\centering
	\includegraphics[width=\textwidth]{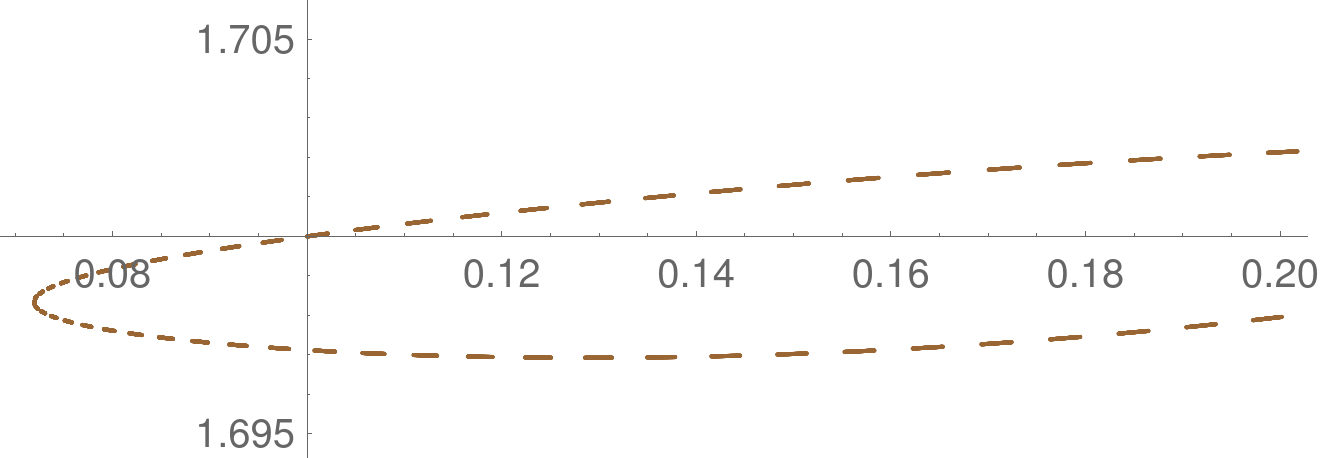}
	\caption{The top brown orbit from Fig.~\ref{fig.orbits18spread}, magnified.}
	\label{fig.TopBrown}
\end{figure}

The top brown orbit in Fig.~\ref{fig.orbits18spread} actually consists of flat broken line ovals; see it magnified in Fig.~\ref{fig.TopBrown} . Its line is broken since 500,000 iterations were insufficient to fill it.

	\begin{figure}[h]
	\centering
	\includegraphics[width=0.7\textwidth]{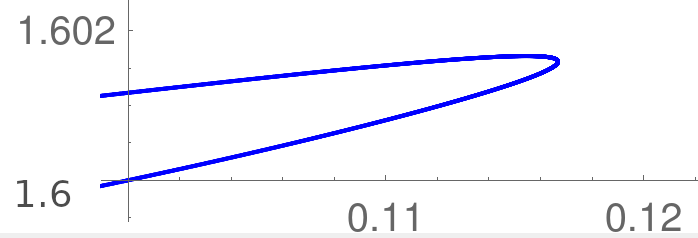}
	\caption{The top blue orbit from Fig.~\ref{fig.orbits18spread}, magnified.}
	\label{fig.TopBlue1}
\end{figure}

The next orbit from the top in Fig.~\ref{fig.orbits18spread} , the blue one, actually consists of very flat ovals; its line is full. See its corner magnified in Fig.~\ref{fig.TopBlue1} ; its ovals are too long and flat to plot one of them fully.

	\begin{figure}[h]
	\centering
	\includegraphics[width=\textwidth]{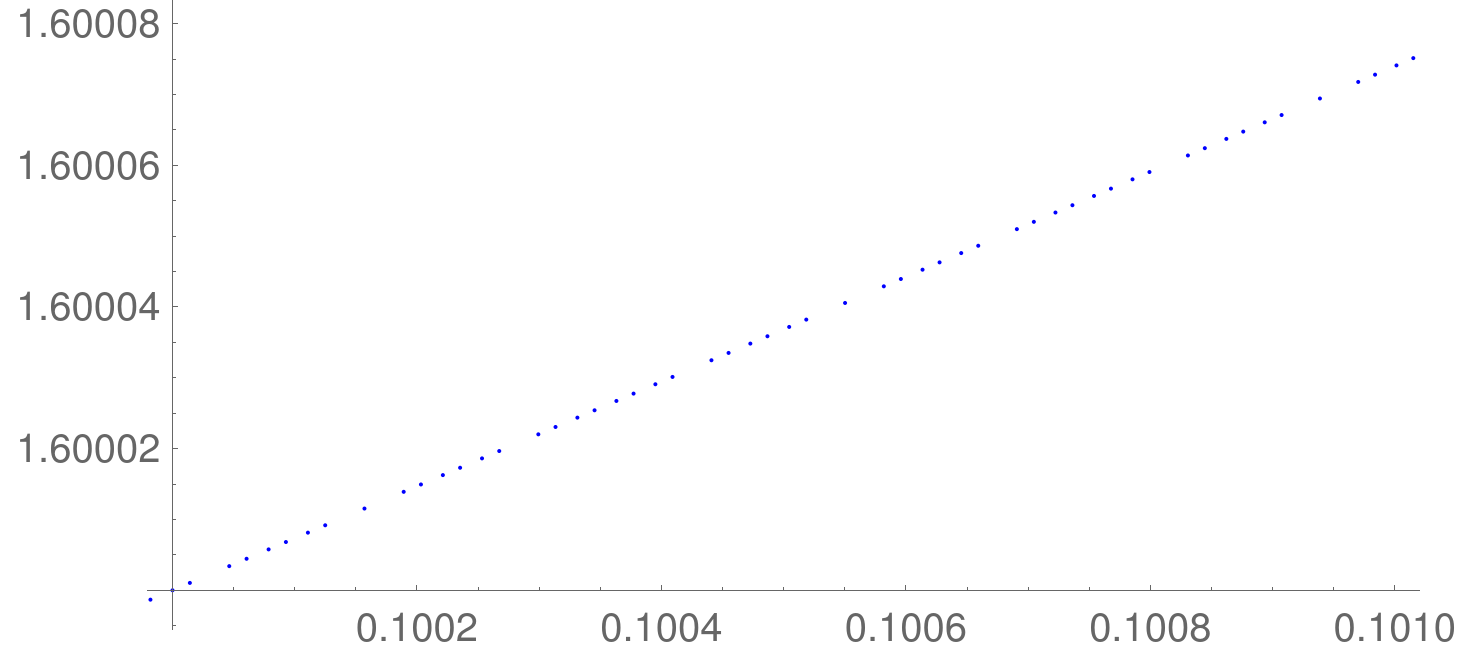}
	\caption{The top blue orbit from Fig.~\ref{fig.orbits18spread}, magnified to individual points.}
	\label{fig.TopBlue2}
\end{figure}

Fig.~\ref{fig.TopBlue2} shows the same orbit magnified to individual points; it shows clearly that these points lie along a curve and the orbit is non-chaotic.  We have plotted all the non-chaotic orbits with similar magnification to see that their points lie along a curve; we do not show all these plots here.

	\begin{figure}[h]
	\centering
	\includegraphics[width=\textwidth]{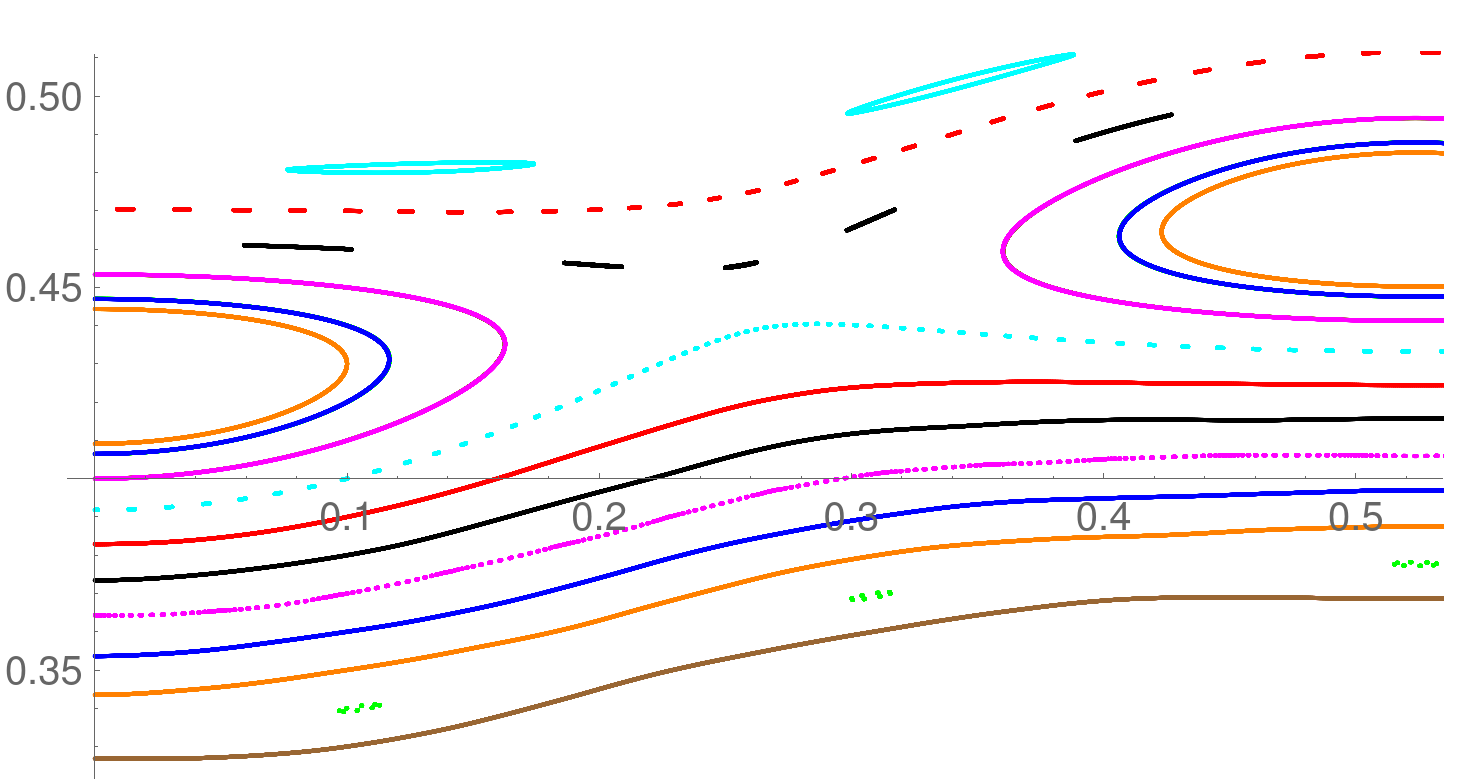}
	\caption{16 orbits spread over the region of Fig.~18 of \cite{fetter}.}
	\label{fig.FetterFig18Orbits}
\end{figure}

In Fig.~\ref{fig.FetterFig18Orbits} we show 16 orbits with the initial points $\phi=0.1$ and $p=0.33,0.34,...,0.48$; these orbits are spread over the region that approximately corresponds to Fig.~18 of \cite{fetter}. These orbits were also calculated through 500,000 iterations with 30 decimals precision.  

In fact, this plot does not show all of these orbits since the orbit with $p=0.45$ (the magenta ovals on the plot) obscures a very close orbit with $p=0.41$, and the orbit with $p=0.44$ (the blue ovals) obscures a very close orbit with $p=0.42$. Pieces of the same color on this plot which are spread over in the horizontal direction, are parts of the same orbit. Each of the dots and dashes on this plot represents an individual oval, except for the green clusters on the bottom, where each of the green clusters is a broken line oval. We have magnified each of the orbits to the level of individual points and found that they lie on one curve in each of the cases, indicating no chaos.

	\begin{figure}[h]
	\centering
	\includegraphics[width=0.6\textwidth]{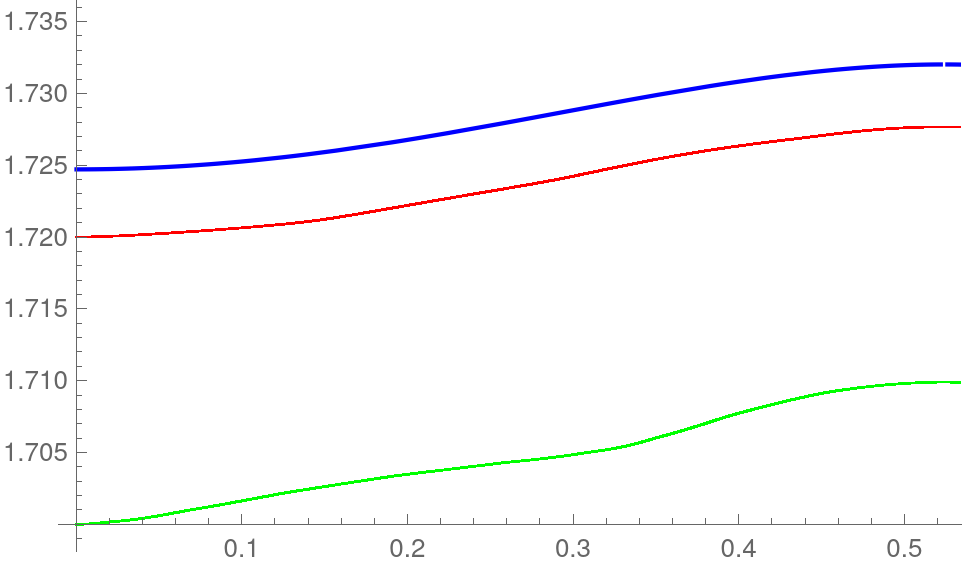}
	\caption{Two orbits, shown in {\color{green}\bf green} and {\color{red}\bf red}, near the boundary of the phase cylinder (shown in {\color{blue}\bf blue}).}
	\label{fig.orbitsNearBdary}
\end{figure}

We took two of the orbits shown in Fig.~\ref{fig.FetterFig18Orbits} , the top red one and the top black one, and calculated them through over 100,000,000 iterations with 80 decimals, with the concurrent control calculation with 100 decimals and continuous monitoring of the discrepancy: if the latter exceeded $10^{-15}$, the calculation was to abort. We hoped to see some slow development of chaos, and these orbits were chosen because their magnified plots gave some mild indications of such possibility.  However, we have not seen any kind of slow chaos and found that the individual points of these orbits lie along smooth curves after all these iterations.

In Fig.~\ref{fig.orbitsNearBdary} we show two orbits near the boundary of the phase cylinder; they are non-chaotic.  The green orbit was calculated through 400,000 iterations while the red one through 3,000,000; both with 30 decimals.  Control calculations with 40 decimals were performed through the first 200,000 iterations for the green orbit (the discrepancy was $\sim 10^{-18}$) and through the first 50,000 iterations for the red orbit (the discrepancy was $6.5\cdot 10^{-20}$).

	\begin{figure}[h]
	\centering
	\includegraphics[width=0.6\textwidth]{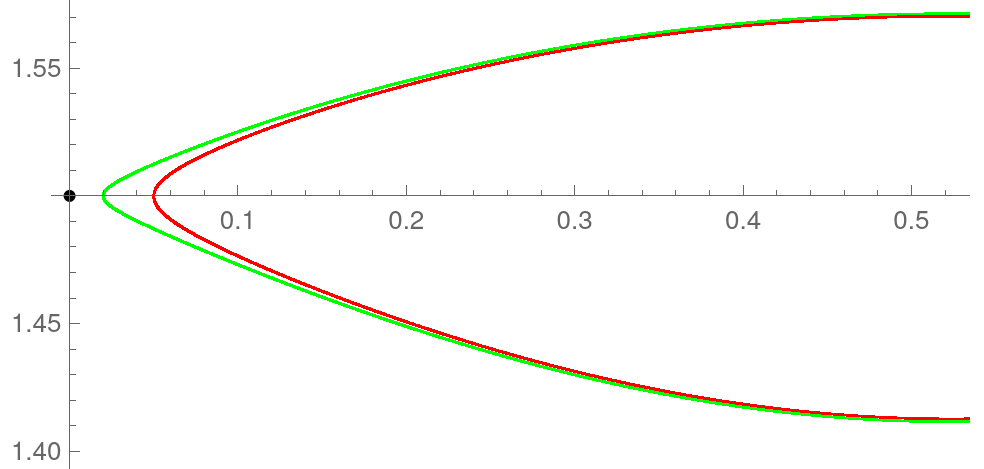}
	\caption{Two orbits, shown in {\color{green}\bf green} and {\color{red}\bf red}, near a hyperbolic periodic point of order 6 ($\phi=0, p=1.5$, shown in {\color{black}\bf black}).}
	\label{fig.TwoOrbitsNearHyperbolic1.5}
\end{figure}

In Fig.~\ref{fig.TwoOrbitsNearHyperbolic1.5} we show two orbits near a hyperbolic periodic point of order 6 ($\phi=0, p=1.5$); they are non-chaotic. The red orbit was calculated through 200,000 iterations while the green one through 3,000,000; both with 30 decimals.  Control calculations with 40 decimals were performed through all 200,000 iterations for the red orbit (the discrepancy was $1.3\cdot 10^{-19}$) and through the first 50,000 iterations for the green orbit (the discrepancy was $3\cdot 10^{-20}$).

\subsection{Chaos near $\frac{5}{12}$}
We have found chaotic orbits near the hyperbolic periodic orbit of period 12 with rotation number $\frac{5}{12}$; one of the points of this orbit has coordinates $(0.2612330773570752, 0.44734020841438477)$.

	\begin{figure}[h]
	\centering
	\includegraphics[width=\textwidth]{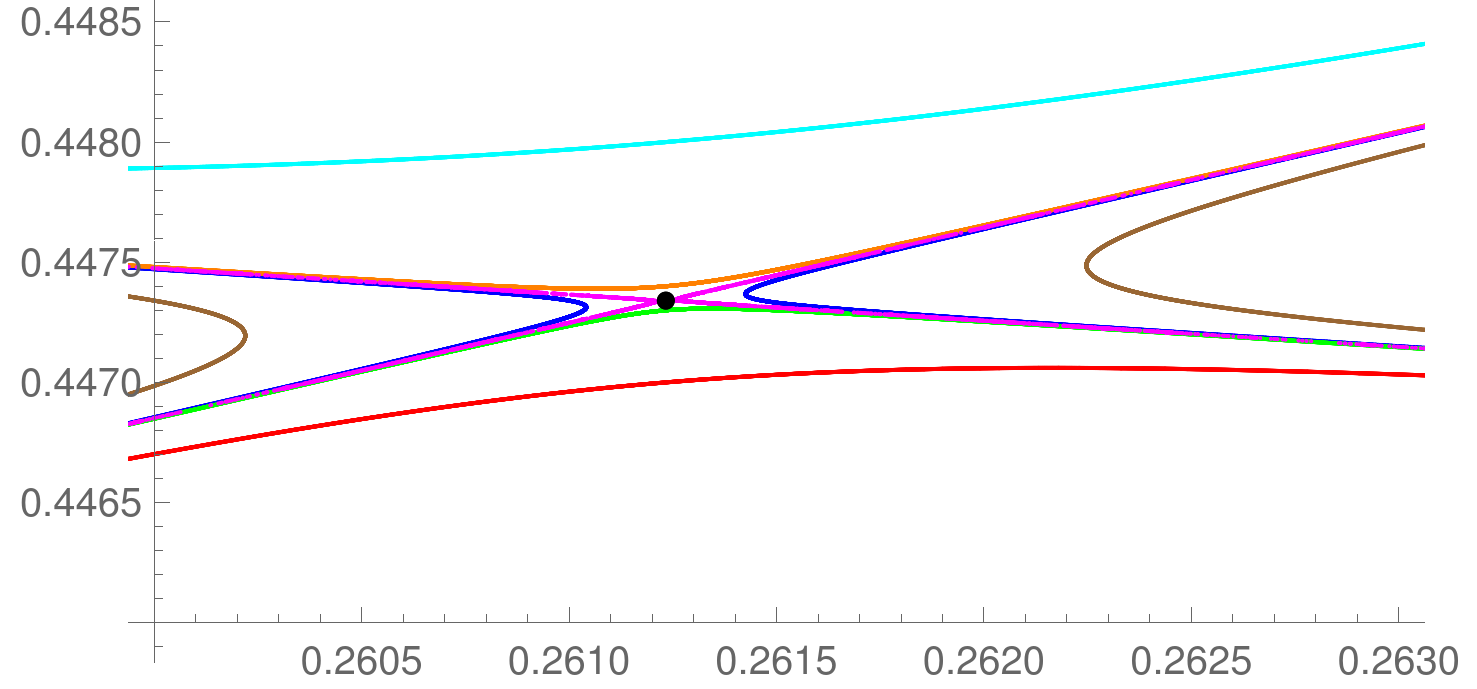}
	\caption{Orbits near the hyperbolic $\frac{5}{12}$ orbit (shown in {\color{black}\bf black}), including one chaotic one.}
	\label{fig.chaos0}
\end{figure}

Fig.~\ref{fig.chaos0} shows orbits near this point (shown as a black dot).  All orbits shown on this figure, have been calculated through 500,000 iterations with 30 decimals, except for the innermost magenta orbit which has been calculated with 50 decimals.  The control calculations have been performed for all these orbits through 50,000 iterations and they have showed the discrepancies of no more than $1.5\cdot 10^{-14}$.

	\begin{figure}[h]
	\centering
	\includegraphics[width=\textwidth]{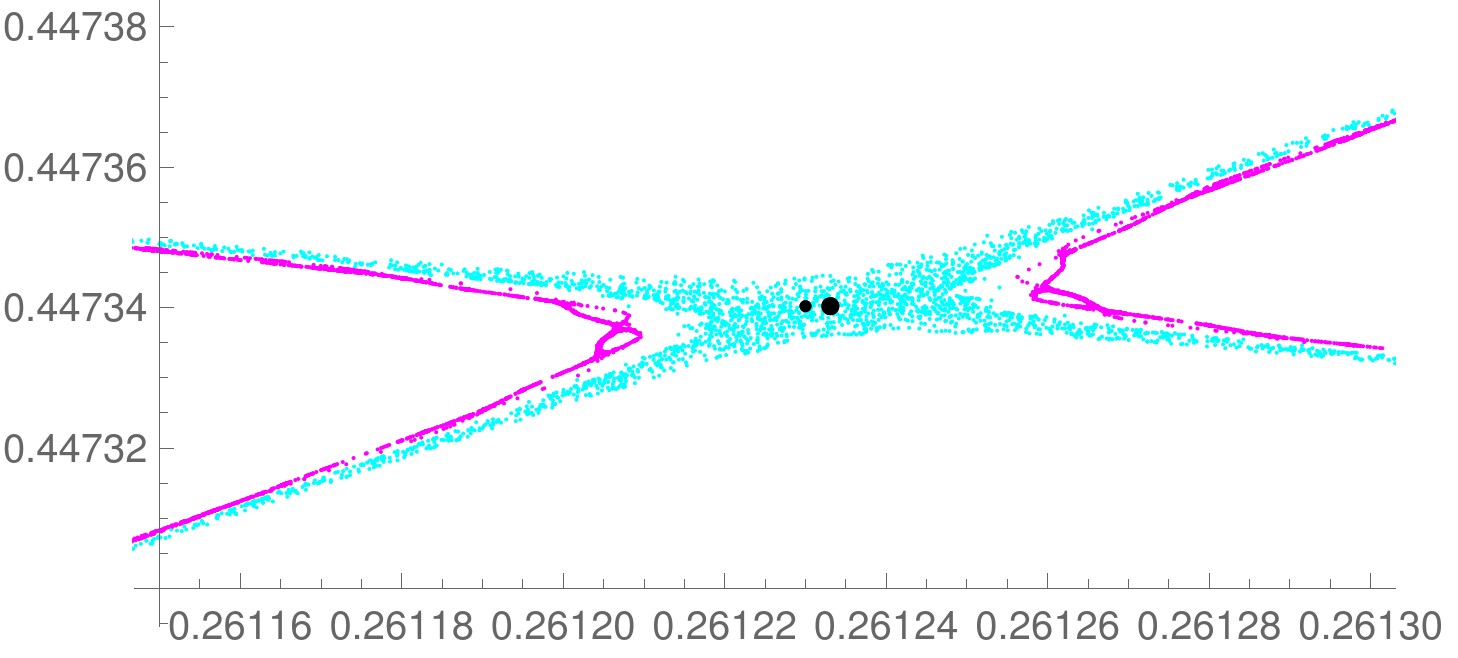}
	\caption{The innermost {\color{magenta}\bf magenta} orbit from Fig.~\ref{fig.chaos0} magnified, together with a definitely chaotic orbit closer to the hyperbolic periodic point (in  {\color{cyan}\bf cyan}). The hyperbolic periodic orbit is shown as a larger {\color{black}\bf black} circle, while the initial point of the cyan orbit is shown as a smaller {\color{black}\bf black} circle.}
	\label{fig.chaos2}
\end{figure}

We have checked that all these orbits are non-chaotic (magnification to the point level shows that the points lie on a curve) except for the innermost magenta orbit which is chaotic, as we see in Fig.~\ref{fig.chaos2} .

\begin{figure}[h]
	\centering
	\includegraphics[width=\textwidth]{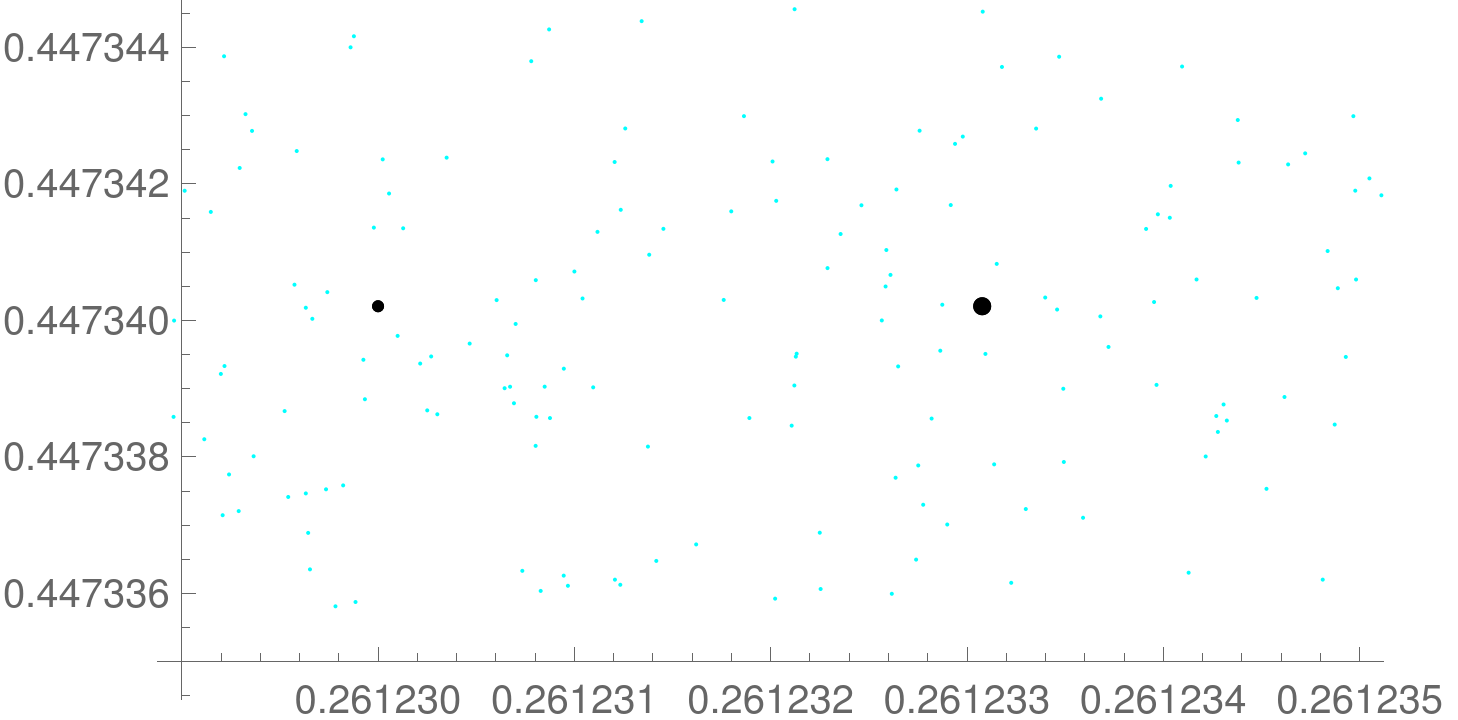}
	\caption{Fig.~\ref{fig.chaos2} magnified further.}
	\label{fig.chaos3}
\end{figure}

We have calculated one more orbit, closer to the hyperbolic periodic orbit, and it shows definite chaos; see Fig.~\ref{fig.chaos2} and~\ref{fig.chaos3} . This orbit has been calculated with the precision of 1,000 decimals, and a control calculation has been performed through the same number of iterations with 1,050 decimals.  The discrepancy has gradually increased with the number of iterations to $10^{-35}$ for 140,000 iterations and 3 for 150,000 iterations; we show 140,000 in Fig.~\ref{fig.chaos2} and~\ref{fig.chaos3} .

\subsection{Chaos near $\frac{11}{27}$}
We have found another chaotic orbit near a hyperbolic periodic orbit of period 27 with the rotation number $\frac{11}{27}$.

This periodic orbit is located between the top cyan ovals in Fig.~\ref{fig.FetterFig18Orbits} ; we found chaos near the point $P=(0.0647611040986, 0.4806888855137)$ on this periodic orbit. (In fact, there are two such hyperbolic periodic orbits near these cyan ovals, one close to their left corners and the other close to their right corners.  These orbits are obtained one from the other by the symmetry with respect to the vertical line $\phi=\pi/6$.  Similarly, applying this symmetry to the cyan ovals yields another orbit consisting of flat ovals located between the cyan ovals. The hyperbolic periodic orbits are located in the remaining small spaces between all these ovals.)

	\begin{figure}[h]
	\centering
	\includegraphics[width=\textwidth]{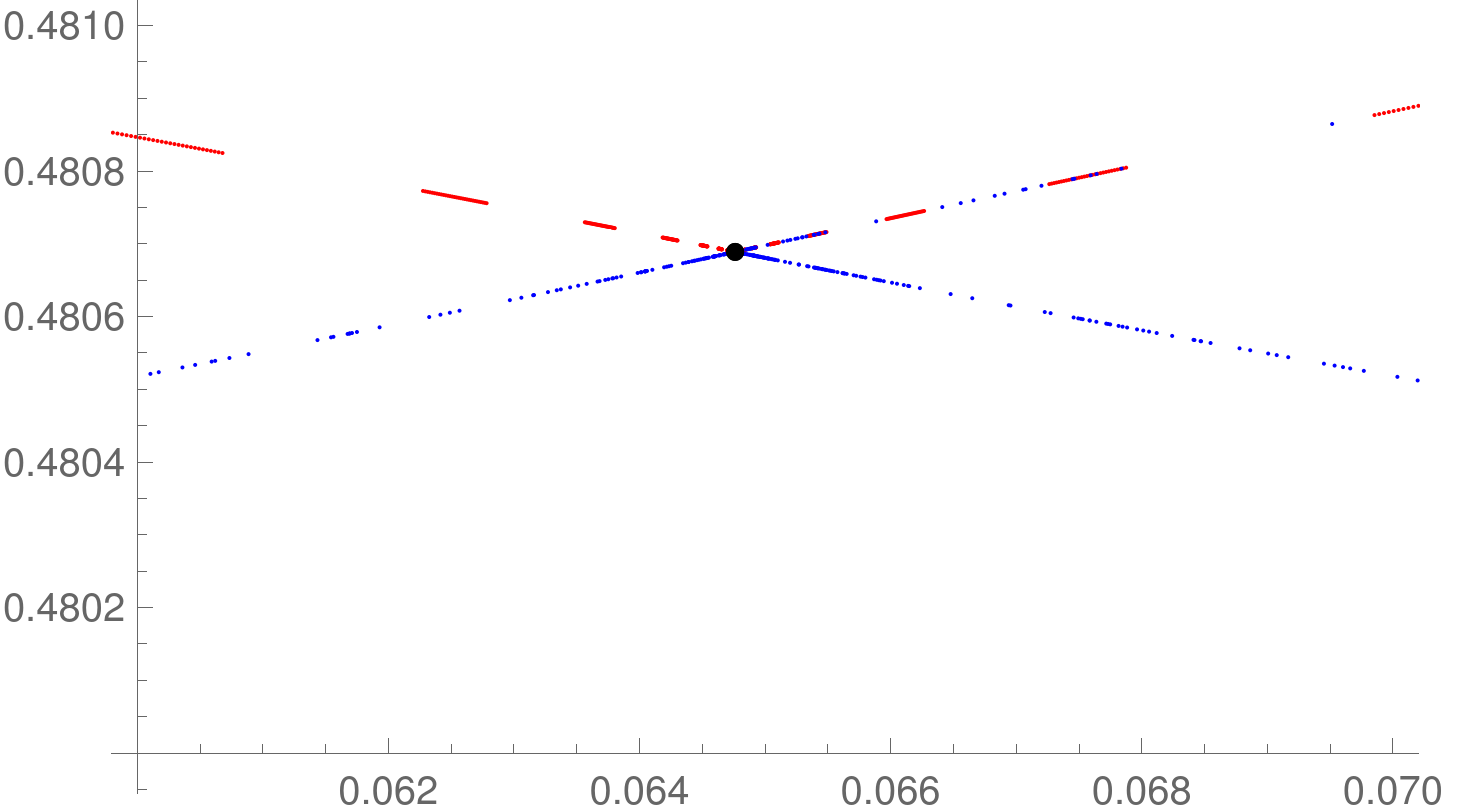}
	\caption{Orbits near the hyperbolic periodic point $\frac{11}{27}$.}
	\label{fig.chaos10}
\end{figure}

In Fig.~\ref{fig.chaos10} we show two orbits near $P$, which is shown in {\color{black}\bf black}:
\begin{description}
\item[{\color{red}\bf red}] The initial point is (0.06476110449368037, 0.480689); it is above the hyperbolic point. The orbit was calculated through 50,000 iterations with 50 decimals.  The discrepancy of the control calculation with 60 decimals was $10^{-27}$.
\item[{\color{blue}\bf blue}] The initial point is (0.06476110449368037, 0.4806888855); it is below the hyperbolic point. The orbit was calculated through 50,000 iterations with 150 decimals. The discrepancy of the control calculation with 170 decimals was $10^{-56}$.
\end{description}

	\begin{figure}[h]
	\centering
	\includegraphics[width=\textwidth]{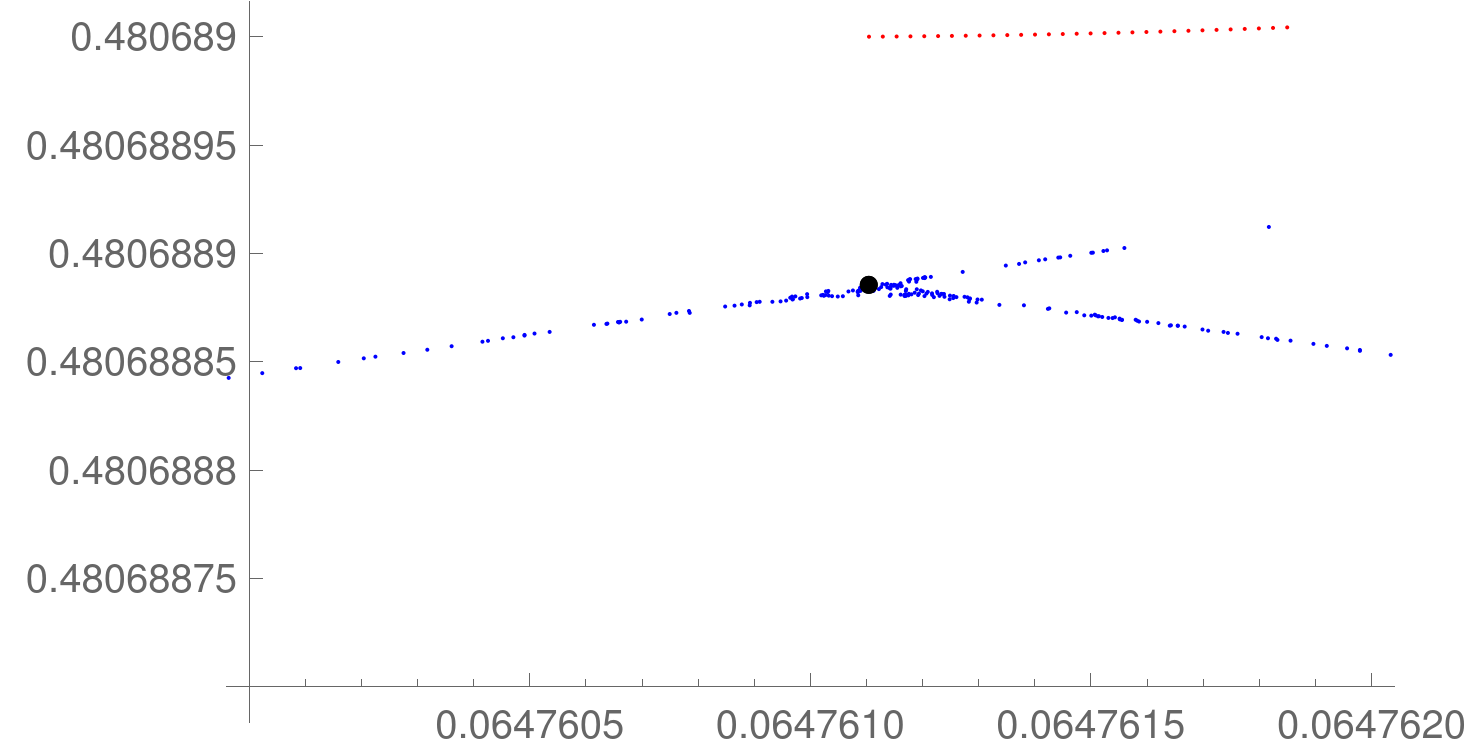}
	\caption{Fig.~\ref{fig.chaos10} magnified}
	\label{fig.chaos11}
\end{figure}

In Fig.~\ref{fig.chaos11} we see that the red orbit is non-chaotic while the blue one shows signs of chaos.

	\begin{figure}[h]
	\centering
	\includegraphics[width=\textwidth]{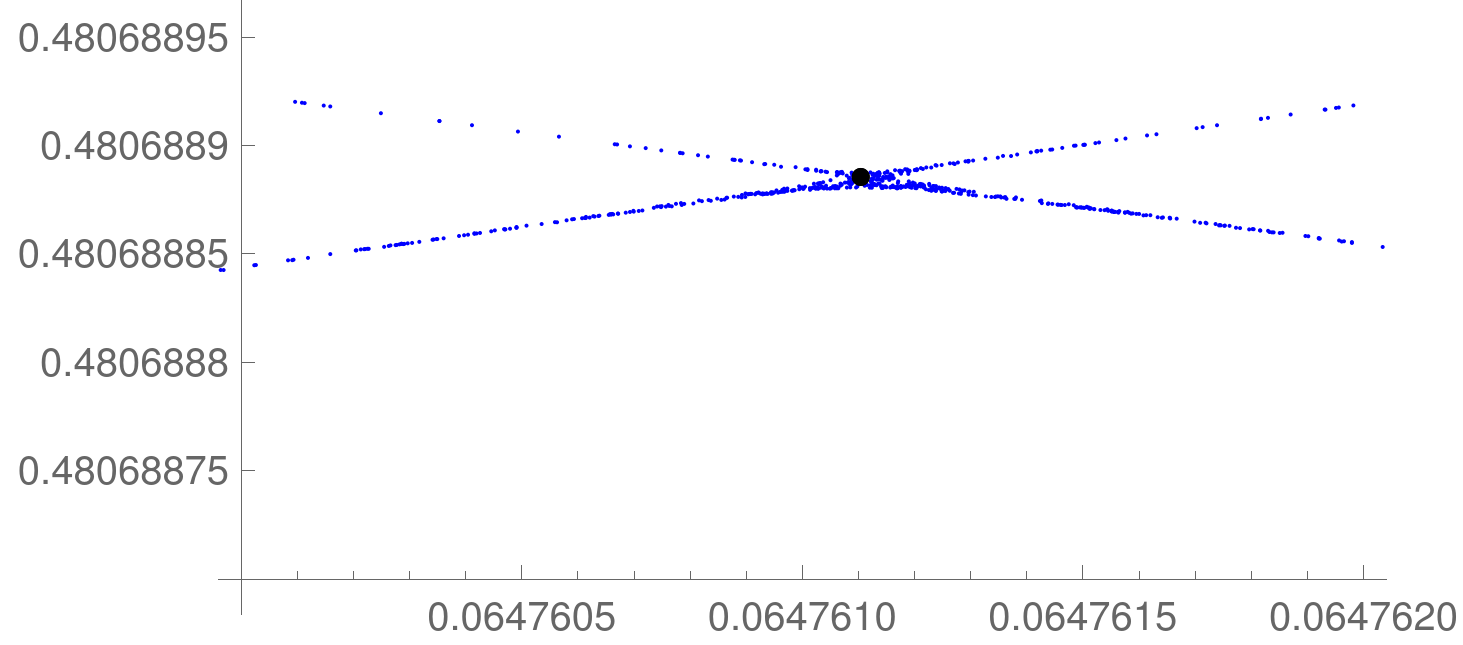}
	\caption{The blue orbit of Fig.~\ref{fig.chaos10} calculated through 150,000 iterations with 500 decimals.}
	\label{fig.chaos12}
\end{figure}

	\begin{figure}[h]
	\centering
	\includegraphics[width=\textwidth]{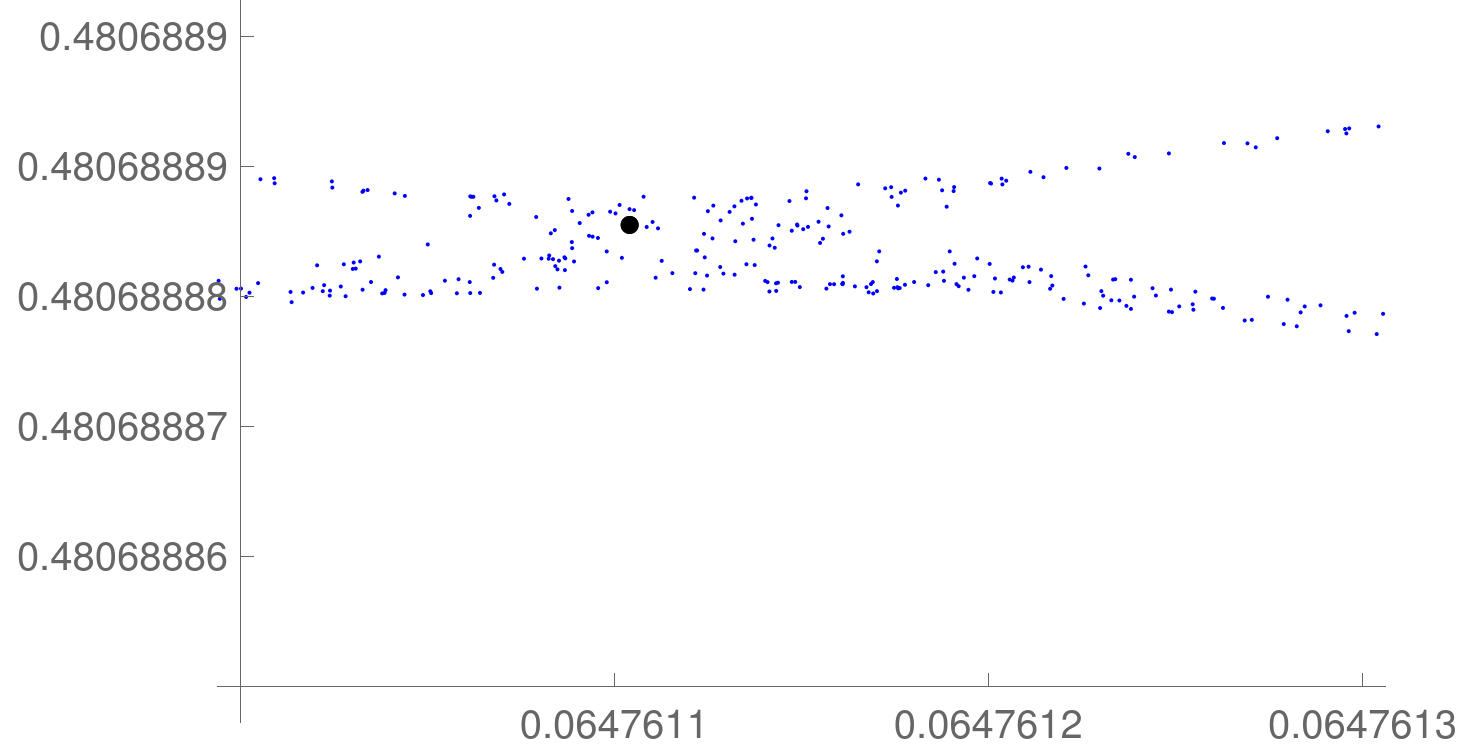}
	\caption{Fig.~\ref{fig.chaos12} magnified}
	\label{fig.chaos13}
\end{figure}

Fig.~\ref{fig.chaos12} and~\ref{fig.chaos13} show the blue orbit calculated through 150,000 iterations with 500 decimals; the discrepancy of the control calculation with 550 decimals was $10^{-234}$.  This orbit is definitely chaotic.


\clearpage
\begin{thebibliography}{90}
	

	\bibitem{AB} M. Arnold, M. Bialy. {\it Nonsmooth convex caustics for Birkhoff billiards.} Pacific J. Math. 295 (2018), no. 2, 257--269. 
	
	\bibitem{BM}M. Bialy, A.E.Mironov, 
	{\it The Birkhoff-Poritsky conjecture for centrally-symmetric billiard tables.}
	Preprint arXiv:2008.03566
	
	%
	
	\bibitem{BT} M. Bialy, S. Tabachnikov, Dan Reznik's identities and more. arXiv:2001.08469
	
	
	

	
\bibitem {fetter}Fetter, Hans L. {\it Numerical exploration of a hexagonal string billiard.} Phys. D 241 (2012), no. 8, 830--846.
	\bibitem{K-S}V. Kaloshin, A. Sorrentino. {\it On the local Birkhoff Conjecture for convex billiards},  Ann. of Math., Vol. 188 (2018), 315--380.
	
	\bibitem{K-S1}V. Kaloshin, A. Sorrentino. {\it On the integrability of Birkhoff billiards.},
	Phil. Trans. R. Soc. A 376: 2017 0419.
	http://dx.doi.org/10.1098/rsta.2017.0419.
	

	
	\bibitem{Poritsky} H. Poritsky. {\it The billiard ball problem on a table with a convex boundary--an illustrative dynamical problem.} Ann. of Math. (2) 51 (1950), 446--470. 
	
\bibitem {Simo} Simó, C.; Vieiro, A. Dynamics in chaotic zones of area preserving maps: close to separatrix and global instability zones. Phys. D 240 (2011), no. 8, 732–753. 
	\bibitem{Tab} S. Tabachnikov. {\it Geometry and billiards.}  Amer. Math. Soc., Providence, RI, 2005.
	
	
	
	

\end{thebibliography}
\end{document}